\documentclass[a4paper,12pt]{amsart}
\usepackage[letterpaper, margin=1.2in]{geometry}
\usepackage{latexsym}
\usepackage{amssymb}
\usepackage{amsmath}
\usepackage{amsfonts}
\usepackage{xcolor}
\usepackage{color}

\newtheorem{thm}{Theorem}[section]

\newtheorem{lem}[thm]{Lemma}

\newtheorem{cor}[thm]{Corollary}

\theoremstyle{definition}

\newtheorem*{rem}{Remark}

\def\fph{\mathbb{F}_{\ph}}

\newcommand{\N}{\mathbb N}
\newcommand{\Z}{\mathbb Z}
\newcommand{\z}{\mathbb Z}
\newcommand{\Q}{\mathbb Q}

\newcommand{\F}{\mathbb F}
\newcommand{\fp}{\mathbb F_p}

\newcommand{\q}{\mathbb Q}
\def\F{\mathbb{F}}

\newcommand{\p}{\mathfrak{p}}

\def\al{\alpha}
\def\la{\lambda}

\def\th{\theta}
\def\md#1{\ \mbox{\rm(mod }{#1})}
\def\nph#1{N_{\ph}(#1)}
\def\npp#1{N_{\ph}^+(#1)}
\def\ph{\phi}

\newcounter{cs}
\stepcounter{cs}
\newcommand{\casos}{\begin{itemize}}
\newcommand{\fcasos}{\end{itemize}\setcounter{cs}{1}}

\newfont{\tit}{cmr12 scaled \magstep3}

\begin{document}
\title[]{On the irreducible  factors of a polynomial
}
\textcolor[rgb]{1.00,0.00,0.00}{}
\author{Lhoussain El Fadil}
\address{Faculty of Sciences Dhar El Mahraz, P.O. Box  1874 Atlas-Fes , Sidi mohamed ben Abdellah University,  Morocco}\email{lhouelfadil2@gmail.com}
\keywords{ Irreducibility criterion,  irreducible factors, Newton polygon techniques} \subjclass[2010]{11R09,
12E05}
\date{\today}
\begin {abstract}
 Jakhar shown that  for $f(x)=a_nx^n + a_{n-1}x^{n-1}+\cdot+ a_0$ ($a_0\neq 0$) is a polynomial with rational coefficients,  if there exists a prime integer $p$ satisfying $\nu_p(a_n)=0$ and $n\nu_p(a_i)\ge (n-i)\nu_p(a_0)> 0$ for every $0\le  i\le  n-1$, then $f(x)$ has at most $gcd(\nu_p(a_0),n)$ irreducible factors over the field $\Q$ of rational numbers and each irreducible factor has degree at least $n/gcd(\nu_p(a_0),n)$. The goal of this paper is to generalize this criterion in the following context: Let $(K,\nu)$ be a rank one discrete valued field, $R_\nu$ its valuation ring and $\F_\nu$ its residue field. {Assume that  $f(x)=\ph^n(x) + a_{n- 1}(x)\ph^{n-1}(x)+\cdot+ a_0(x)\in R_\nu[x]$, with for every $i=0,\dots,n-1$, $a_i(x)\in R_\nu[x]$,  and $a_0(x)\neq 0$ for some monic polynomial  $\ph\in R_\nu[x]$ with $\overline{\ph}$ is irreducible in $\F_\nu[x]$. If for every $0\le  i\le  n-1$,  $n\nu_p(a_i)\ge (n-i)\nu_p(a_0)>0$,} then $f(x)$ has at most $gcd(\nu_p(a_0(x)),n)$ irreducible factors over the field $K^h$ and so over $K$ and each irreducible factor has degree at least $n/gcd(\nu_p(a_0),n)$, where $K^h$ is the henselization of $(K,\nu)$. 
 \end{abstract}
\maketitle
 \section{Introduction}
 Polynomial factorization over a field is very useful in algebraic number theory (prime ideal factorization), coding theory (generators of cyclic codes), Galois theory, extension of valuations, etc. Irreducible polynomials appear naturally in the study of polynomial factorization and the question of irreducibility of polynomials over henselian fields  has been of interest to many mathematicians (cf. \cite{Ei, Gi, Mu, We, H, O}). In 1850, Eisenstein gave the most popular irreducibility criterion \cite{Ei}. In 1906, Dumas gave a generalization to this criterion  as follows if 
 $f(x)=a_nx^n + a_{n-1}x^{n-1}+\dots+ a_0$, with $a_0\neq 0$ is a polynomial with rational coefficients and if there exists a prime integer $p$ satisfying $\nu_p(a_n)=0$ and $n\nu_p(a_i)\ge (n-i)\nu_p(a_0)> 0$ for every $0\le  i\le  n-1$ and $gcd(\nu_p(a_0),n)=1$, then $f(x)$ is irreducible over $\Q$ \cite{Ei}.
 In 2008, R. Brown \cite{Br} gave a simple proof of the most general version of Eisenstein-Scho\"nemann irreducibility Criterion. Namely, if $f(x)=\ph^n(x) + a_{n -1}(x)\ph^{n-1}(x)+\cdot+ a_0(x)\in R_\nu[x]$, with for every $i=0,\dots,n-1$, $a_i(x)\in R_\nu[x]$, deg$(a_i)<$deg$(\ph)$, and gcd$(\nu(a_0(x),n)=1$ for some monic polynomial  $\ph\in R_\nu[x]$ with $\overline{\ph}$ is irreducible in $\F_\nu[x]$, then $f(x)$ is irreducible  over the field $K$.   In 2020, Jakhar relaxes the condition  gcd$(\nu(a_0,n)=1$, in the context of  $\ph=x$, and he shown that for a polynomial $f(x)=a_nx^n + a_{n-1}x^{n-1}+\dots+ a_0$,   with rational coefficients and  $a_0\neq 0$, if there exists a prime integer $p$ satisfying $\nu_p(a_n)=0$ and $n\nu_p(a_i)\ge (n-i)\nu_p(a_0)> 0$ for every $0\le  i\le  n-1$, then $f(x)$ has at most $gcd(\nu_p(a_0,n)$ irreducible factors over the field  over $\Q$ (see \cite{Ja}).
In this paper, we extend both Anuj's and Eisenstein-Scho\"nemann's irreducibility criterion \cite{Br, Ja}, as follows that : Let $(K,\nu)$ be a rank one discrete valued field, $R_\nu$ its valuation ring, $M_\nu$ its maximal ideal, and $\F_\nu$ its residue field. {Assume that  $f(x)=\ph^n(x) + a_{n- 1}(x)\ph^{n-1}(x)+\cdot+ a_0(x)\in R_\nu[x]$, with for every $i=0,\dots,n-1$, $a_i(x)\in R_\nu[x]$,  and $a_0(x)\neq 0$ for some monic polynomial  $\ph\in R_\nu[x]$ with $\overline{\ph}$ is irreducible in $\F_\nu[x]$. If for every $0\le  i\le  n-1$,  $n\nu_p(a_i)\ge (n-i)\nu_p(a_0)>0$}, then $f(x)$ has at most $gcd(\nu_p(a_0(x)),n)$ irreducible factors over the field $K^h$. In particular, let $f(x)\in \Z[x]$ be an irreducible polynomial over $\Q$ and $L=\Q(\al)$ be the number field generated by a complex root $\al$ of $f(x)$. Assume that there exists  a prime integer $p$,  with $\overline{f(x)}=\overline{\ph(x)}^n$ in $\F_p[x]$ for some monic polynomial $\ph$, whose reduction is irreducible over $\F_p$. Let  $f(x)=\ph^n(x) + a_{n-1}(x)\ph^{n-1}(x)+\cdot+ a_0(x)$ be the $\ph$-expansion of $f(x)$ satisfying  for every $i=0,\dots,n-1$, $\nu_p(a_i)\ge \nu_p(a_0)-i\la$, with $\la=\frac{\nu_p(a_0)}{n}$, then $f(x)$ has at most $gcd(\nu_p(a_0(x)),n)$ irreducible factors over the field $\Q_p$. In particular, there are at most $gcd(\nu_p(a_0(x)),n)$ prime ideal of $\Z_K$ lying above $p$.  Our results are illustrated by examples.
\section{Notations}\label{not}
As our proof is based on Newton polygon's techniques, we recall some fundamental notations and techniques  on Newton polygon.
 For any number field $L=\Q(\al)$ generated by a complex root $\al$ of a monic irreducible polynomial $f(x)\in \Z[x]$, in 1894, K. Hensel
developed a powerful approach by showing that the prime ideals of ${\mathbb Z}_L$  
lying above a prime $p$ are in one-to-one correspondence with
irreducible factors of $f(x)$ in ${\mathbb Q}_p[x]$, where ${\mathbb Z}_L$ is the ring of integers of $L$. For every prime ideal
corresponding to any irreducible factor in ${\mathbb Q}_p[x]$, the
ramification index  and the residue degree together are the same as
those of the local field defined  by the irreducible factor
\cite{H}.  This result was generalized in \cite[Proposition 8.2]{Ne} as follows for a valued field  $(K,\nu)$ and $L=K(\al)$ a simple algebraic extension generated  by $\al\in \overline{K}$ a root of a monic irreducible polynomial $f(x)\in R_\nu[x]$, the valuation of $L$ extending $\nu$  are in one to one correspondence with irreducible factors of $f(x)$ in $K^h[x]$.
So, in order to describe all valuations of $L$ extending $\nu$, we need to factorize the polynomial $f(x)$ into irreducible factors in $K^h[x]$. The first step of the factorization is based on Hensel's lemma. Unfortunately, the factors provided by Hensel's lemma are not necessarily irreducible over $K^h$. The Newton polygon techniques could refine the factorization. Namely, the theorem of the product, the theorem of the polygon and the theorem of the residual polynomial say that we can factorize any factor provided by Hensel's lemma, with as many sides of the polygon and as many irreducible factors of each residual polynomial. For more details, see \cite{GMN, O} for Newton polygons over $p$-adic numbers and \cite{CMS, El}  for Newton polygons over rank one discrete valuations.  
 For any  monic polynomial 
$\phi\in R_\nu[x]$ {\it whose reduction} modulo $M_\nu$ is irreducible  in
${\mathbb F}_\nu[x]$, let ${\mathbb F}_{\phi}$ be 
the field $\frac{{\mathbb F}_\nu[x]}{(\overline{\phi})}$.

Let $(K,\nu)$ be a rank one discrete valued field, $R_\nu$ its valuation ring, $M_\nu$ its maximal ideal, $\F_\nu$ its residue field,  and $(K^h,\nu^h)$ its henselization.  By normalization, we can assume that $\nu(K^*)=\Z$ and so $M_\nu$  is a principal ideal of  $R_\nu$ generated by an element $\pi\in K$ satisfying $\nu(\pi)=1$. Let $f(x)\in R_\nu[x]$ be a monic polynomial and assume that
$\overline{f(x)}$ is a power of $\overline{\ph}$ in $\F_\nu[x]$, with $\overline{\ph}$ is irreducible in $\F_\nu[x]$.
 Upon to  the Euclidean division
 by successive powers of $\phi$, we can expand $f(x)$ as follows $f(x)=\sum_{i=0}^la_i(x)\phi(x)^{i}$, called    the $\phi$-expansion of $f(x)$
 (for every $i$, deg$(a_i(x))<$
deg$(\phi)$). 
The $\phi$-Newton polygon of $f$, denoted by $\nph{f}$ is lower boundary convex envelop of the set of points $\{(i,\nu(a_i)),\, i=0,\dots,n\}$ in the Euclidean plane.  For every edge $S_j$, of the polygon, let  $A_{j-1}=(i_{j-1}, \nu(a_{i_{j-1}}))$ its initial point   and $A_j=(i_{j}, \nu(a_{i_{j}}))$ its final point. 
 Let $l(S_j)=i_{j}-i_{j-1}$ be its length,  $h(S_j)=\lambda_jl(S_j)$  its height, and  $\lambda_j=\frac{\nu_p(a_{i_{j}})-\nu_p(a_{i_{j-1}})}{i_{j}-i_{j-1}}\in {\mathbb Q}$,  called the slope of $S_j$.  
 Remark that $\nph{f}$ is the process of joining the obtained edges $S_1,\dots,S_r$ ordered by  the increasing slopes, which  can be expressed as $\nph{f}=S_1+\dots + S_r$.  {The segments  $S_1,\dots,$ and $S_r$ are called the sides of $\nph{f}$.} 
 For every side $S$ of the polygon $\nph{f}$, $l(S)$ is 
 the length of its projection to the $x$-axis and  $h(S)$  is the length of its projection to the $y$-axis.

 For every side $S$ of $\nph{f}$, with initial point $(s, u_s)$ and length $l$, let $d=l/e$, called the degree of $S$. For every 
$0\le i\le l$, we attach   the following
{{\ residual coefficient} $c_i\in{\mathbb F}_{\phi}$:
$$c_{i}=
\left
\{\begin{array}{ll} 0,& \mbox{ if } (s+i,{\it u_{s+i}}) \mbox{ lies strictly
above } S
 \mbox{ or }
{\it u_{s+i}}=\infty,\\
\left(\dfrac{a_{s+i}(x)}{p^{{\it u_{s+i}}}}\right)
\,\,
\mbox{\rm(mod }{(\pi,\phi(x))}),&\mbox{ if }(s+i,{\it u_{s+i}}) \mbox{ lies on }S.
\end{array}
\right.$$
where $(\pi,\phi(x))$ is the maximal ideal of $R_\nu[x]$ generated by $\pi$ and $\phi$. That means if $(s+i,{\it u_{s+i}}) \mbox{ lies on }S$, then $c_i=\overline{\dfrac{a_{s+i}(\beta)}{\pi^{{\it u_{s+i}}}}}$, where $\beta$ is a root of $\phi$.}
 
Let $\lambda=-h/e$ be the slope of $S$, where  $h$ and $e$ are positive coprime integers, and let $d=l/e$ be the degree of $S$.  Notice that, 
the points  with integer coordinates lying in $S$ are exactly $(s,u_s),(s+e,u_{s}+h),\cdots, (s+de,u_{s}+dh)$. Thus, if $i$ is not a multiple of $e$, then 
$(s+i, u_{s+i})$ does not lie in $S$, and so, $c_i=0$. Let
$f_S(y)=t_0y^d+t_1y^{d-1}+\cdots+t_{d-1}y+t_{d}\in{\mathbb F}_{\phi}[y]$ be 
the residual polynomial of $f(x)$ associated to the side $S$, where for every $i=0,\dots,d$,
 $t_i=c_{ie}$.
{
\begin{rem}
 Note that if $\nu(a_{s+i}(x))=0$ and $\phi=x$, then ${\mathbb F}_{\phi}={\mathbb F}_\nu$ and $c_i=\overline{{a_{s+i}}} \mbox{\rm(mod }{\pi})$. Thus this notion of residual coefficient generalizes the reduction modulo a maximal ideal. If $\lambda=0$, then for every $i=0,\dots,d$, $(s+i,{\it u_{s+i}}) \mbox{ lies on }S$ if and only if $\nu(a_{s+i}(x))=0$. Thus if $\lambda=0$ and $\phi=x$, then $c_i=\overline{{a_{s+i}}} \mbox{\rm(mod }{\pi})$ and $f_S(y)\in{\mathbb F}_\nu[y]$ coincides with the reduction of $f(x)$ modulo the maximal ideal $M_\nu=(\pi)$.
    \end{rem}}
The  theorem of the product and theorem of the polygon play a key role in the proof (see \cite[Lemma 3.1, 3.3, and Th. 3.5]{El}).
\section{Main results}
Let $(K,\nu)$ be a rank one discrete valued field, $R_\nu$ its valuation ring, $M_\nu$ its maximal ideal, $\F_\nu$ its residue field,  and $(K^h,\nu^h)$ its henselization.  Let $f(x)\in R_\nu[x]$ be a monic polynomial and assume that
$\overline{f(x)}$ is a power of $\overline{\ph}$ in $\F_\nu[x]$, with $\overline{\ph}$ is irreducible in $\F_\nu[x]$.
 
 According to notations and terminologies of section \ref{not}, Eisenstein-Scho\"nemann irreducibility Criterion and  Jakhar's criterion could reformulated as follows:\\
{\bf Eisenstein-Scho\"nemann irreducibility Criterion:}\\ 
 {\it {Let $f(x)\in \Z[x]$ be a polynomial. If for some prime integer $p$, $\overline{f(x)}=\overline{\ph}^l$ for some monic polynomial $\ph\in\Z[x]$, whose reduction is irreducible over $\fp$, and  $\nph{f}=S$ has a single side of degree $d=1$, with respect to $\ph$ and $\nu=\nu_p$, then $f(x)$ is irreducible over  $\Q$.}}\\
 {\bf Jakhar's criterion:}\\ 
 {\it {Let $f(x)\in \Z[x]$ be a polynomial. If for some prime integer $p$, $\nph{f}=S$ has a single side, with respect to $\ph=x$ and $\nu=\nu_p$, then $f(x)$  has at most $d$ monic irreducible factors in $\Q[x]$, where $d$ is the degree of $S$.}}\\
 Eisenstein-Scho\"nemann irreducibility Criterion could be generalized as follows:
 \begin{thm}\label{irred}(\cite[Cor. 3.2]{El})\\
 Let $f(x)\in R_\nu[x]$ be a monic polynomial such that
$\overline{f(x)}$ is a power of $\overline{\ph}$ in $\F_\nu[x]$, with $\overline{\ph}$ is irreducible in $\F_\nu[x]$.
Let  $f(x)=\ph^n(x) + a_{n-1}(x)\ph^{n-1}(x)+\cdot+ a_0(x)$ be the $\ph$-expansion of $f(x)$. If for every $i=0,\dots,n-1$,  $\nu(a_i)\ge \nu(a_0)-i\la$, with $\la=\frac{\nu(a_0)}{n}$ and $f_S(y)$ is irreducible over $\fph$, then $f(x)$ is irreducible  over $K^h$.
  \end{thm}
Jakhar's irreducibility criterion could be generalized as follows (it relaxes $\ph=x$ required in Jakhar's criterion).
\begin{thm}\label{bound}
 Let $f(x)\in R_\nu[x]$ be a monic polynomial such that
$\overline{f(x)}$ is a power of $\overline{\ph}$ in $\F_\nu[x]$, with $\overline{\ph}$ is irreducible in $\F_\nu[x]$.
Let  $f(x)=\ph^n(x) + a_{n-1}(x)\ph^{n-1}(x)+\cdot+ a_0(x)$ be the $\ph$-expansion of $f(x)$. If for every $i=0,\dots,n-1$,  $\nu(a_i)\ge \nu(a_0)-i\la$, with $\la=\frac{\nu(a_0)}{n}$, then $f(x)$ has at most gcd$(\nu(a_0(x)),n)$ irreducible factors over the field $K^h$  of degree at least { $em$, where $e=n/d$ and deg$(\ph)=m$.} 
\end{thm}
\begin{cor}\label{val}
Under the hypothesis of Theorem \ref{bound}, assume that $f(x)$ is irreducible over $K$ and let $L=K(\al)$ be the simple extension generated by $\al\in \overline{K}$ a root of $f(x)$. If for every $i=0,\dots,n-1$,  $\nu(a_i)\ge \nu(a_0)-i\la$, with $\la=\frac{\nu(a_0)}{n}$, then $L$ has at most $d=$gcd$(\nu(a_0(x)),n)$ distinct valuation of $L$ extending $\nu$.
\end{cor}
\begin{cor}\label{prime}
Under the hypothesis of theorem \ref{bound}, assume that $f(x)\in \Z[x]$ is irreducible over $\Q$ and let $L=\Q(\al)$ be the number field  generated by $\al$ a complex root of $f(x)$ and $\overline{f(x)}$ is a power of $\overline{\ph}$ in $\F_p[x]$ for some  prime integer $p$ and a monic polynomial $\ph\in\Z[x]$, whose reduction $\overline{\ph}$ is irreducible in $\F_p[x]$. If for every $i=0,\dots,n-1$,  $\nu(a_i)\ge \nu(a_0)-i\la$, with $\la=\frac{\nu(a_0)}{n}$, then $\Z_L$ has at most $gcd(\nu_p(a_0(x)),n)$ distinct prime ideals   lying above  $p$, where $\Z_L$ is the ring of integers of  $L$.
\end{cor}
\begin{thm}\label{bound2}
 Let $f(x)\in R_\nu[x]$ be a monic polynomial and let
$\overline{f(x)}=\prod_{i=1}^r\overline{\ph_i}^{n_i}(x)$ be the factorization  of $\overline{f(x)}$ in $\F_\nu[x]$, with every $\ph_i\in R_\nu[x]$ is a monic polynomial. For  every $i=1,\dots, r$, let $N_{\ph_i}(f)=S_{i1}+\dots+S_{ig_i}$ be the $\ph_i$-Newton polygon of $f(x)$.  For  every $i=1,\dots, r$ and  $j=1,\dots, g_i$, let  $d_{ij}=\frac{l_{ij}}{e_{ij}}$ be the degree of $S_{ij}$, where $l_{ij}$ is the length of $S_{ij}$, $-\frac{h_{ij}}{e_{ij}}$ is its slope, $h_{ij}$ and $e_{ij}$ are two positive coprime integers.
Then $f(x)$ has at most $\sum_{i=1}^r\sum_{j=1}^{g_i}d_{ij}$ irreducible factors over the field $K^h$ and so, over $K$.
\end{thm}
\begin{cor}\label{val2}
Under the hypothesis of Theorem \ref{bound2}, assume that $f(x)$ is irreducible over $K$ and let $L=K(\al)$, where $\al\in \overline{K}$ is a root of $f(x)$. Then there at most $\sum_{i=1}^r\sum_{j=1}^{g_i}d_{ij}$ distinct valuations of $L$ extending $\nu$.
\end{cor}
\begin{cor}
Let $R$ be a Dedekind domain with quotient field $K$ and $\p$ be a nonzero prime ideal of $R$. Let $\nu$ be the $\p$-adic valuation of $K$  and 
$\overline{f(x)}=\prod_{i=1}^r\overline{\ph_i}^{n_i}(x)$ be the factorization  of $\overline{f(x)}$ in $\F_\nu[x]$, with every $\ph_i\in R_\nu[x]$ is a monic polynomial.  Under the hypothesis of  Corollary \ref{bound2}, assume that $f(x)$ is irreducible over $K$ and let $L=K(\al)$, where $\al\in \overline{K}$ is a root of $f(x)$. Then there at most $\sum_{i=1}^r\sum_{j=1}^{g_i}d_{ij}$ prime ideals of $\Z_L$ lying above $\p$, where $\Z_L$ is the integral closure of $R$ in $L$.
\end{cor}
\section{Proofs}
\begin{proof} of of Theorem \ref{bound}.
Under the hypothesis of Theorem \ref{bound}, let $f(x)=f_1(x)\times\dots\times f_t(x)$ be the factorization of $f(x)$ in $K^h[x]$, with every $f_i(x)$ is monic.
By Gauss's lemma, every $f_i(x)\in R_\nu[x]$.  Let $i=1,\dots,t$. Since $\overline{f_i(x)}$ divides   $\overline{f(x)}$, then $\overline{f_i(x)}=\overline{\ph^{l_i}(x)}$ for some natural integer $l_i$. The hypothesis of  Theorem \ref{bound}, implies that $\nph{f}=S$ has a single side of slope $-\la$. By the theorem of the product \cite{El}, for every $i=1,\dots,t$, $\nph{f_i}=S_i$ has a single side of slope $-\la$, 
  $S=S_1+\dots+S_t$, and $f_S(y)=\prod_{i=1}^t{f_i}_{S_i}(y)$ up to multiply by a non zero element of $\fph$.  Let $\la=h/e$,  where $h$ and $e$ are two positive coprime integers. Since the length $l(S)=ed$ and $l(S)=n=\sum_{i=1}^tl(S_i)=e\times(d_1+\dots+d_t)$, then $d=d_1+\dots+d_t$. The fact that $d_i\ge 1$ for every $i=1,\dots,t$, implies that $d\le t$.  {Fix $i=1,\dots,t$.  As $\nph{f_i}=S_i$, then  $f_i(x)=\ph^{l_i}+\dots+b_{i0}$ is the $\ph$-expansion of $f_i(x)$, with $l_i=d_ie$. Thus deg$({f_i}(x))=med_i=me_i$ deg$({f_i}_{S_i}(y))$. As deg$({f_i}_{S_i}(y))=d_i\ge 1$,  deg$({f_i}(x))\ge me$ as desired.}
\end{proof}
\begin{proof} of Theorem \ref{bound2}.
Under the hypothesis of Theorem \ref{bound2}, by Gauss's lemma and  Hensel's lemma, $f(x)=f_1(x)\times\dots\times f_r(x)$ in $R_{\nu^h}[x]$ such that for every $i=1,\dots,r$,   $\overline{f_i(x)}=\overline{\ph_i^{l_i}(x)}$. Fix $i=1,\dots,r$. 
By theorem of the product and theorem of the polygon \cite{El}, $f_i(x)=\prod_{j=1}^{g_i}f_{ij}(x)$ such that for every $j=1,\dots,g_i$, $N_{\ph_i}(f_{ij})=S_{ij}$ has  single side of slope $-\la_{ij}=-h_{ij}/e_{ij}$, where $h_{ij}$ and $e_{ij}$ are two positive coprime integers.
By Theorem \ref{bound}, every  $f_{ij}$  has at most  $d_{ij}$ monic irreducible factor in $K^h[x]$, and finally we get the desired result.
\end{proof}
\section{Examples}
\begin{enumerate}
\item
Let $f(x)=\ph^6+24x\ph^4+24\ph^3+15(16x+32)\ph+48$ and $\ph\in Z[x]$ be a monic polynomial whose reduction is irreducible in $\F_2[x]$. For $p=2$, $\nph{f}=S$ has a single side of length $l=6$, height $H=4$, and so $d=2$. By Theorem \ref{bound}, $f(x)$ has  at most $2$ irreducible factors in $\q_2[x]$.\\
As $f_S(y)=y^2+y+1$ is irreducible over $\F_2[x]$, then for $\ph=x+m$, with $m\in\Z$,   $f_S(y)$ is irreducible over $\fph\simeq \F_2$. Thus Theorem \ref{irred}, $f(x)$ is irreducible over $\Q_2$. Let $L=\Q(\al)$ and $\Z_L$ its ring of integers, where $\al$is a complex root of $f(x)$. Since $f(x)$ is irreducible over $\Q_2$, there is a single prime ideal of $\Z_L$ lying above $2$.
\item
Let $f(x)=\ph^6+24x\ph^3+9(16x+32)\ph+3(16x+16)$ and $\ph=x^2+x+1\in\Z[x]$. For $p=2$, $\overline{\ph}$ is irreducible over $\F_2$, $\nph{f}=S$ has a single side of length $l=6$, height $H=4$, and  $d=2$. By Theorem \ref{bound}, $f(x)$ has  at most $2$ irreducible factors in $\q_2[x]$.
\item
Let $f(x)=(\ph_1^5+p^3)(\ph_2^4+p^3)$, with for every $i=1,2$, $\ph_i$ is irreducible over $\F_p$, and $\overline{\ph_1}\neq \overline{\ph_2}$.  Then  $\overline{f(x)}=\overline{\ph_1}^5\overline{\ph_2}^4$ in $\F_p[x]$. By combining  Hensel's lemma and Theorem \ref{bound2}, $f(x)$ has exactly two distinct monic  irreducible factors in $\Q_p[x]$.  
 \end{enumerate}
 
\end{document}